\newcommand{\be}{\begin{eqnarray}}
\newcommand{\ee}{\end{eqnarray}}

\def\ll#1{\left#1}
\def\r#1{\right#1}
\def\fr{\frac{1}{2}}
\def\mref#1{(\ref{#1})}
\def\mb{\mbox{}}

\def\p{\partial}
\def\bd{\begin{displaymath}}
\def\ed{\end{displaymath}}
\def\ba#1{\begin{array}{#1}}
\def\ea{\end{array}}
\def\nn{\nonumber}
\newfont{\Bbb}{msbm10 scaled 1200}
\documentclass[12pt]{article}
\usepackage{fleqn}

\title{ Deformed Galilei symmetry
}
\author{
     Piotr Kosi\'nski$^*$,                  Pawe\l\
Ma\'slanka\thanks{supported by KBN grant 2P03B 130 12}\\
                              Theoretical Physics Department II\\
                                     University of \L\'od\'z\\
                          ul. Pomorska 149/153, 90 236 \L\'od\'z,
Poland}
\date{}
\begin{document}
\maketitle
\thispagestyle{empty}
\begin{abstract}
A particular deformation of central extended Galilei group is
considered.
It is shown that the deformation influences the rules of constructing
the
composed systems while one particle states remain basically unaffected.
In
particular the mass appeared to be non additive.
\end{abstract}
\section{Introduction}
Quantum groups have emerged in physics in connection with an a attempt
to
understand the symmetries underlying the exact solvability of certain
quantum--mechanical
and statistical models~\cite{b1}.

In this respect they appeared to be quite powerful. It is therefore
natural to pose the question whether quantum groups could provide a
suitable tool
for describing another symmetries in physics, in particular the
space--time
ones. From the mathematical point of view one can define a variety of
deformations
of classical space--time symmetry groups. They posses some attractive
features.
For example the deformation of Poincare symmetry, constructed
in~\cite{b2}
provides a ten--parameter symmetry ``group" with dimensional parameter
naturally
built in; one can speculate that this parameter is a natural cut--off or
reflects
the change of properties of space--time at very small scale.
However once we take seriously the very idea of quantum symmetries we
are faced
with many conceptual problems: what is the meaning of non--commuting
group 
parameters, how to deal with non--cocommmutative coproduct on the Lie
algebra
level, etc. There exists no general, commonly acceptable solution to
these
questions. An order to shad some light on them in the present paper we
describe
a particular example of rather mild deformation of Galilei group which
admits
a nice and simple physical interpretation (on the Lie algebra level this
deformation was given in~\cite{b3}). The quantum group we are
considering is
a~deformation of the standard central extension of classical Galilei
group
constructed in such a~way that only the relations concerning the
additional
group variable are modified. Our basic assumption is that the projective 
representations of (classical) Galilei group arise from vector
representations
of our deformed group. It appears then that the one--particle
representations
remain unmodified except that the physical mass cannot be directry
identified with
the eigenvalue of mass operator and, moreover is bounded from above. The 
quantum structure influences the way the many--particles states are 
constructed. The main difference is the way the total mass of composed
system is defined; the mass is no longer additive. We address also to
the
problem of non--commutativity of coproduct  for generators and show that
both
direct and transpose coproduct lied to physically equivalent theories.
We hope
that the results presented below show that, at least in cases it is
possible
to give a coherent physical interpretation to the notion of quantum
space--time
symmetry.

\section{The deformed centrally extended\protect\newline
Ga\-li\-lei~group~G$_k$}
Our starting point is the $\kappa$--Poincare group defined in
Ref.\cite{b2}.
Its defining relations read
\be\label{w1}
\mb[{\Lambda^\mu}_\nu,{\Lambda^\alpha}_\beta]&=&0\nn\\
\mb[{\Lambda^\mu}_\nu,a^\rho]&=&\frac{-i}{\kappa}\ll(({\Lambda^\mu}_0-
{\delta^\mu}_0){\Lambda^\rho}_\nu+({\Lambda^0}_\nu-{\delta^0}_\nu)
g^{\mu\rho}\r)\nn\\
\mb[a^\mu,a^\nu]&=&\frac{i}{\kappa}({\delta^\mu}_0a^\nu-{\delta^\nu}_0a^\mu)\\
\Delta({\Lambda^\mu}_\nu)&=&{\Lambda^\mu}_\alpha
\otimes{\Lambda^\alpha}_\nu\nn\\
\Delta(a^\mu)&=&{\Lambda^\mu}_\alpha\otimes a^\alpha+a^\mu\otimes I\nn\\
({\Lambda^\mu}_\nu)^*&=&{\Lambda^\mu}_\nu,\;\;\;\;\;(a^\mu)^*=a^\mu\nn
\ee

In order to obtain the deformed central extension of Galilei group,
G$_k$, we follow
the same procedure as in the classical case. First, we define a trivial
extension
of $\kappa$--Poincare by adding new unitary generator $\xi$,
\be\label{w2}
\xi^*\xi=\xi\xi^*=I\nn\\
\mb[\xi,{\Lambda^\mu}_\nu]=0,\;\;\;[\xi,a^\mu]=0\\
\Delta(\xi)=\xi\otimes\xi\nn
\ee
Then we redefine
\be
\zeta&=&\xi e^{-imca^0}\nn
\ee
$\zeta$\ is again unitary and obeys
\be\label{w3}
a^0\zeta&=&\zeta  a^0\nn\\
a^k\zeta&=&e^{-\frac{mc}{\kappa}}\zeta a^k\nn\\
{\Lambda^0}_0\zeta&=&\zeta\frac{sh(\frac{mc}{\kappa})+{\Lambda^0}_0
ch(\frac{mc}{\kappa})}{ch(\frac{mc}{\kappa})+
{\Lambda^0}_0sh(\frac{mc}{\kappa})}\\
{\Lambda^0}_i\zeta&=&\zeta\frac{{\Lambda^0}_i}{ch(\frac{mc}{\kappa})+
{\Lambda^0}_0sh(\frac{mc}{\kappa})}\nn\\
{\Lambda^i}_0\zeta&=&\zeta\frac{{\Lambda^i}_0}{ch(\frac{mc}{\kappa})+
{\Lambda^0}_0sh(\frac{mc}{\kappa})}\nn\\
{\Lambda^i}_j\zeta&=&\zeta\ll({\Lambda^i}_j-\frac{{\Lambda^i}_0{\Lambda^0}_j}
{1-({\Lambda^0}_0)^2}(cth(\frac{mc}{\kappa}+\alpha)-{\Lambda^0}_0)\r),
\;\;\;cth\alpha\equiv{\Lambda^0}_0\nn
\ee
In order to apply contraction we redefine~\cite{b4}
\begin{eqnarray}\label{w4}
{\Lambda^0}_0&=&\gamma\equiv(1-\vec{v}^2/c^2)^{-\fr}\nn\\
{\Lambda^0}_i&=&\frac{\gamma}{c}v^k{R^k}_i\nn\\
{\Lambda^i}_0&=&\frac{\gamma}{c}v^i\\
{\Lambda^i}_j&=&\ll({\delta^i}_k+(\gamma-1)\frac{v^iv^k}{\vec{v}^2}\r)
{R^k}_j\nn\\
a^0&=&c\tau,\;\;\;a^i\equiv a^i\nn
\end{eqnarray}
where ${R^i}_j$\ represent rotations, $RR^T=R^TR=I$. Then we put
$\kappa\to\infty$,
$c\to\infty$, $\frac{\kappa}{c}\equiv k$--fixed. The resulting structure
reads
\be\label{w5a}
\mb[v^i,v^j]&=&[v^i,{R^k}_j]=[v^i,a^j]=[v^i,\tau]=0\nn\\
\mb[a^i,{R^k}_j]&=&[a^i,a^j]=[a^i,\tau]=0\nn\\
\mb[{R^i}_j,{R^k}_l]&=&[{r^i}_j,\tau]=0\\
\Delta(v^i)&=&{R^i}_j\otimes v^j+v^i\otimes I\nn\\
\Delta(a^i)&=&{R^i}_j\otimes a^j+v^i\otimes\tau+a^i\otimes I\nn\\
\Delta({R^i}_j)&=&{R^i}_k\otimes {R^k}_j\nn
\ee
as well as
\be\label{w5b}
\tau\zeta&=&\zeta\tau\nn\\
a^k\zeta&=&e^{-\frac{m}{k}}\zeta a^k\nn\\
v^k\zeta&=&e^{-\frac{m}{k}}\zeta v^k\\
{R^i}_j\zeta&=&\zeta{R^i}_j\nn\\
\Delta(\zeta)&=&(\zeta\otimes\zeta)e^{-im(\frac{\vec{v}^2}{2}
\otimes\tau+v^k{R^k}_i\otimes a^i)}\nn
\ee

It follows from eq.\mref{w5a} that, as long as the additional generator
$\zeta$\ 
is neglected, we are dealing here with the \underline{classical} Galilei
group.
Therefore we expected that the one--particle states are described by a
standard
theory. However, the description of many--particle states will be
certainly affected.

In order to find, via duality, the corresponding Lie algebra we define
\newline
$\Phi(\zeta,\tau,\vec{a},\vec{v},\vec{\Theta})$, where $\Theta^i$\
parametrise
the rotation matrices R, to be a function on group manifold under the
proviso
that $\zeta$\ stands leftmost. The Lie algebra generators are then
defined by the
following duality rules:
\be\label{w6}
\mb<H,\Phi(\zeta,\tau,\vec{a},\vec{v},\vec{\Theta})>&=&
i\frac{\p\Phi}{\p\tau}|_e\nn\\
\mb<P_k,\Phi(\zeta,\tau,\vec{a},\vec{v},\vec{\Theta})>&=&
-i\frac{\p\Phi}{\p a^k}|_e\nn\\
\mb<J_k,\Phi(\zeta,\tau,\vec{a},\vec{v},\vec{\Theta})>&=&
-i\frac{\p\Phi}{\p\Theta^k}|_e\\
\mb<K_k,\Phi(\zeta,\tau,\vec{a},\vec{v},\vec{\Theta})>&=&
-i\frac{\p\Phi}{\p v^k}|_e\nn\\
\mb<M,\Phi(\zeta,\tau,\vec{a},\vec{v},\vec{\Theta})>&=&\mu
\zeta\frac{\p\Phi}{\p\zeta}|_e\nn
\ee
where $\mu$\ is an arbitrary mass unit (see below).

The standard Hopf algebra duality rules imply then
\be\label{w7a}
\mb[J_i,P_k]&=&i\epsilon_{ikl}P_l,\;\;[J_i,K_k]=i\epsilon_{ikl}K_l\nn\\
\mb[K_i,P_k]&=&i\delta_{ik}\frac{\mu}{1-e^{-2\mu/k}}(1-e^{-\frac{2M}{k}})\\
\mb[K_i,H]&=&iP_i,\nn
\ee
the remaining commutators being vanishing and
\be\label{w7b}
\Delta P_i&=&P_i\otimes a^{-M/K}+I\otimes P_i\\
\Delta K_i&=&K_i\otimes e^{-M/k}+I\otimes K_i\nn
\ee
the remaining coproducts being primitive.

This Hopf algebra, modulo simple redefinitions, was ob\-tained firstly
in Ref.\cite{b3}.
\section{One--particle dynamics}
Due to the fact that the classical Galilei group is a Hopf subalgebra
one can
define in a standard way its projective representations~\cite{b5}
\be\label{w8}
\Phi_\sigma(\vec{p})\to e^{i(-\frac{\vec{p}^2}{2m_f}\tau+\vec{p}
\cdot\vec{a})}\sum\limits_{\sigma'}D^s_{\sigma\sigma'}(R)
\Phi_{\sigma'}(R^{-1}\vec{p}-m_f\vec{v});
\ee
here $m_f$\ denotes the physical mass of a particle. In the classical
case the projective
representations can be converted to the vector representation of a
central extension
of Galilei group~\cite{b5}. Usually, the mass parameter enters the
definition
of this extension so that it depends on the projective representation we
have
selected. However, it is not difficult to rephrase the whole
construction in
such a way that there is a universal central extension which produces
all
projective representations. In fact, we choose an arbitrary reference
mass $\mu$\ to define the central extension and consider the following
representation
\be\label{w9}
\Phi_\sigma(\vec{p})\to \zeta^{\frac{m_f}{\mu}}
e^{i(-\frac{\vec{p}^2}{2m_f}\tau+\vec{p}\cdot\vec{a})}
\sum\limits_{\sigma'}D^s_{\sigma\sigma'}(R)\Phi_{s'}(R^{-1}\vec{p}-m_f\vec{v})
\ee
It is obvious that one obtains the projective representation~\mref{w8}.

We follow the same strategy in the deformed case. Namely, we demand that
{\protect\samepage
\be\label{w10}
\rho\Phi_s(\vec{p})&=&(I\otimes\zeta^{\frac{m}{\mu}})
e^{i(-\frac{\vec{p}}{2m_f}\otimes\tau+\vec{p}\otimes\vec{a})}\nn\\
&&\sum\limits_{\sigma'}(I\otimes D^s_{\sigma\sigma'}(R))
\Phi_{\sigma'}(\vec{p}\otimes R^{-1}-m_f I\otimes\vec{v})
\ee}
be the representation of G$_k$. It follows from the duality rules that
$m$\ is
an eigenvalue of $M$. However, contrary to the classical case, it
differs in 
general from $m_f$. In fact, it is easy to check that $\rho$, as given
by eq.\mref{w10},
provides a representation of G$_k$,
\be\label{w11}
(\rho\otimes I)\circ\rho&=&(I\otimes\Delta)\circ \rho
\ee
only if the following relation holds
\be\label{w12}
m_f&=&\frac{\mu}{1-e^{-2\mu/k}}(1-e^{-2m/k})
\ee
In the limit $k\to\infty$\ we get $m_f=m$, as expected.

At first it seems that the condition that we are dealing with
the representations of G$_k$\ does not give anything new: eq.\mref{w12}
is merely a parametrization of the physical mass $m_f$\ in terms
of some parameter $m$. However, it is $m$\ which is an eigenvalue of
primitive
generator $M$\ so it is $m$\ and not $m_f$\  which is additive.

It is easy to see that by a dimensionless (and $m_f$--independent!)
redefinition of the $K_i$\ and $P_i$\ generators eq.\mref{w12} can be
put in form
\be\label{w13}
m_f&=&\frac{k}{2}(1-e^{-2m/k})
\ee
It follows now that the addition formula for physical masses reads
\be\label{w14}
M_f&=&m_f+m'_f-\frac{2m_fm'_f}{k}
\ee
Eqs.\mref{w13} and \mref{w14} have interesting properties. Eq.\mref{w13}
implies
$m_f<\frac{k}{2}$\ so that $\frac{k}{2}$\ is a maximal possible mass
value. This is the
only modification for one--particle theory.
It is then easy to check that this result agrees with the addition
formula:
$m_f<\frac{k}{2},\;m'_f<\frac{k}{2}$\ imply $M_f<\frac{k}{2}$; moreover,
if $m_f=\frac{k}{2}$\ then $M_f=\frac{k}{2}$; the value $\frac{k}{2}$\
plays the 
role of infinite mass in the theory (see also below).
\section{Two--particle dynamics}
In order to describe two--particle states we demand that they should
also transform according to the representations of G$_k$. Let us first
remind briefly
the relevant construction in the undeformed case. The Galilean
generators are then simply the 
sums of generators corresponding to both particles, $P_i=P_{1i}+P_{2i}$,
etc. 
Applying this assumption to the energy operators one obtain $H=H_1+H_2$\ 
which implies that the particles do not interact. However, due to the
specific form of Galilei Lie algebra this situation can be cured in a
simple way (contrary to the relativistic case). The energy operator
does not appear on the right--hand side of commutation rules. Therefore,
one can redefine $H=H_1+H_2\to H_1+H_2+V$, where $V$\ commutes with all
generators. Due to the fact that two--particle representation is
reducible, $V$\ 
can be nontrivial. In fact, it can be any function of ``relative"
dynamical
variables
$(\vec{p}_1-\vec{p}_2)^2,\;(\frac{\vec{K}_{1}}{m_{f_1}}-\frac{\vec{K}_2}{m_{f_2}})^2$,
etc.;
these relative variables form, together with Galilean generators, the
complete
set of dynamical variables.

Let us follow the same procedure in the deformed case. Using the
coproduct
formulae~\mref{w7a} we define the generators of two--particle
representation
of G$_k$:
\be\label{w15}
\vec{P}&=&e^{-\frac{m'}{k}}\vec{p}+\vec{p}'\\
\vec{K}&=&e^{-{m' \over k}}\vec{k}+\vec{k}'\nn
\ee
the remaining generators being simply the sums of their one--particle
counterparts.
Let us put for a moment aside the problem concerning noncommutativity of
the 
coproduct (the construction of total dynamical variables seems to depend
on numbering of particles). In order to construct the two--particle
dynamics we 
supply the generators of G$_k$\ by additional ones (relative variables)
such
that (i) together with the generators of G$_k$ they form a complete set
of
operators, (ii) they have definite properties under rotations and (iii)
the standard Heisenberg commutation rules among coordinate and momenta
hold.

It is easy to check that these conditions are met by the following
choice
\be\label{w16}
\vec{p}&=&e^{-m'/k}\vec{p}+\vec{p}\kern2pt'\;\;\;\;\;\mbox{(\rm total
momentum})\nn\\
\vec{R}&=&{1 \over M_f}(e^{-m'/k}\vec{K}+\vec{K}')\;\;\;\;\;\mbox{\rm
(``center--of--mass" coordinate)}\\
\vec{\Pi}&=&{1\over M_f}(m_f'\vec{p}-m_f
e^{-m'/k}\vec{p}\kern2pt')\;\;\;\;\;\mbox{\rm (relative momentum)}\nn\\
\vec{\rho}&=&{\vec{K}\over m_f}-e^{-m'/k}{\vec{K}'\over
m_f'}\;\;\;\;\;\mbox{\rm (relative coordinate)}\nn
\ee
Let us note the following identity
\be
H&=&{\vec{p}\kern2pt^2\over
2m_f}+{\vec{p\kern1pt\lower3pt\hbox{$'$}}^2\over 2m_f'}={\vec{P}'\over
2M_f}+{\vec{\Pi}^2\over 2v_f}\nn
\ee
where
\be\label{w17}
v_f&=&{m_fm_f'\over M_f}
\ee
is the deformed reduced mass. Let us note that $v_f=m_f$\ if
$m_f'=\frac{k}{2}$\ 
which supports the point of view that $k\over 2$\ plays a role of
infinite mass.
Now, following the standard procedure we can introduce the interaction
by adding
to the total kinetic energy on arbitrary function
$V(|\vec{\rho}|,|\vec{\Pi}|,\vec{\rho}\cdot\vec{\Pi})$.
If, as in the undeformed case we restrict ourselves to the
coordinate--dependent
potential functions we arrive at the following form of total hamiltonian
\be\label{w18}
H&=&{\vec{P}^2\over 2M_f}+\ll({\vec{\Pi}^2\over 2v_f}+V(|\vec{\rho}|)\r)
\ee

As an illustration consider a simple toy model of the k--deformed
hydrogen atom.
Its hamiltonian reads:
\be\label{w19}
H&=&{\vec{p}\kern1.5pt^2\over
2m_f}+{\vec{p\kern1pt\lower4pt\hbox{$'$}}^2\over 2m_f'}-{e^2\over
\ll|\vec{r}-\sqrt{1-{2m_f'\over k}}\vec{r}\kern1pt'\r|}
\ee
where $m_f$\ and $m_f'$\ are the masses of electron and proton,
respectively
(their role can be exchanged -- see below). It follows from
eq.\mref{w18}
that in the relative coordinates $H$\ takes a standard form
\be\label{w20}
H&=&{\vec{p}\kern1.5pt^2\over 2M_f}+{\vec{\Pi}^2\over 2v_f}
-{e^2\over |\vec{\rho}|}
\ee
so the energy spectrum reads
\be\label{w21}
E_n&=&-{v_fe^4\over 2\hbar^2n^2}=-{ve^4\over 2\hbar^2n^2}
\ll(1+{2v\over k}+\ldots\r)
\ee
where $v$\ is standard reduced mass $\ll(={m_fm_f'\over m_f+m_f'}\r)$.
The 
question can be posed whether this correction is in principle
observable. This
depends on whether $m_f$\ and $m_f'$\ are observable masses of electron
resp.~proton.
They should be measured by making them interacting with infinitely heavy
source 
of forces which, in our theory means that its mass equals $k\over 2$.
But in this
case the reduced mass is just the mass of the particle under
consideration.
Therefore, the masses $m_f$\ and $m_f'$\ are, in principle, measurable
and
so is the correction to the Bohr formula.

Let us now consider the problem of apparent asymmetry in construction of
two--particle
states (the problem of coproduct in the context of particle interaction
was
also considered in~\cite{b6}). It seems that our theory depends on the
choice
of order in which we add particles to the system;  this is obviously
related
to the noncocommutativity of the coproduct. If it were true to the whole
theory
would make no sense. However, it is not difficult to show that the order
in
which we add particles is immaterial -- both description are related by
an
unitary transformation. In fact, using the transposed coproduct we
arrive at
the following set of basic dynamical variables replacing those given by
eqs.\mref{w16}
\be\label{w22}
\tilde{\vec{p}}&=&\vec{p}+e^{-m/k}\vec{p}\kern1.5pt'\nn\\
\tilde{\vec{R}}&=&{1\over M_f}\ll(\vec{K}+e^{-m/k}\vec{K}\kern1pt'\r)\\
\tilde{\vec{\Pi}}&=&{1\over M_f}\ll(m_f'\vec{p}e^{-m/k}-m_f
\vec{p}\kern1.5pt'\r)\nn\\
\tilde{\vec{\rho}}&=&{\vec{K}e^{-m/k}\over m_f}
-{\vec{K}\kern1pt'\over m_f'}\nn
\ee
First two equations represent simply the transposed coproduct. The
relative
variables are chosen  such that they obey Heisenberg commutation rules.
The
ambiguity in sign of $\tilde{\vec{\Pi}}$\ and $\tilde{\vec{\rho}}$\ is
resolved
by demanding that our transformation does not exchange particles (i.e.
reduces
to an identity for $k\to\infty$). It is easy to check that both sets of
variables,
eqs.\mref{w16} and~\mref{w22} are related by the following unitary
transformation
\be\label{w23}
U&=&e^{{i\over\sqrt{m_fm_f'}}\arctan\ll({\sqrt{m_fm_f'}
(1-\sqrt{1-\frac{m_f}{k}}\sqrt{1-\frac{m_f'}{k}})\over
m_f\sqrt{1-\frac{2m_f'}
{k}}+m_f'\sqrt{1-\frac{2m_f}{k}}}\r)(\vec{K}\otimes\vec{P}-
\vec{P}\otimes\vec{K})}
\ee
We conclude that whatever order of composing the system we choose, the
result is,
up to an unitary transformation, the same. Let us now consider a
particular case of identical 
particles. Define the exchange operator $S$\ as
\be\label{w24}
S\Phi_{\sigma\sigma'}(\vec{p},\vec{p}\kern1.5pt')&=&
\Phi_{\sigma'\sigma}(\vec{p}\kern1.5pt',\vec{p})
\ee
It easy to see that
\be\label{w25}
(US)^2&=&I
\ee
Therefore, one can define bosonic vs. fermonic states as obeying
\be\label{w26}
(US\Phi)_{\sigma\sigma'}(\vec{p},\vec{p}\kern1.5pt')&=&
\pm\Phi_{\sigma\sigma'}(\vec{p},\vec{p}\kern1.5pt')
\ee
In the limit $k\to\infty,\;\;\;U\to I$\ and our definition coincides
with
the standard one. For any $k$\ it reverses the sign of relative
coordinates.
Consequently, it imposes the same restriction on admissible states and
observables as classical Fermi--Bose symmetry.
\section{Conclusions}
We have presented simple nonrelativistic quantum mechanical model based
on deformed
centrally extended Galilei group. The deformation considered Leaves the
proper
Galilei group unaffected. The main assumption was that the physical
states form
the vector representations of this central extension. The one--particle
states
do not differ from their undeformed counterparts and span the
representation
space of projective representation of standard Galilei group, the only
difference
being the relation between the physical mass and the eigenvalue of
additional
generator (mass operator). In particular, this relation implies an upper
bound
for particle mass, $m_f<{k\over 2}$. The main difference concerns the
way the 
many--particle states are constructed. The addition formula for masses
is 
modified in the way it contains the deformation parameter $k$. It obeys
a consistency condition that adding two masses below the ${k\over
2}$--bound
one gets the total mass obeying the same inequality. Also, the
definition of
reduced mass is modified. However, we have shown that once the
two--particle
interaction is expressed in terms of relative coordinates (it should be 
expressible if the theory is Galilei--invariant) the theory takes
standard form except the above --- mentioned difference in the
definition
of reduced mass. Corrections to the reduced mass are, in principle,
observable
because the masses of individual particles can be determined by making
them interacting with ``infinitely heavy" ($m_f={k\over 2}$) source.

The main difficulty inherent in the formalism is that the total as well
as
relative variables depend on the order in which particles are added to
the
system. However, the theories obtained by selecting different orders are
unitary equivalent which makes the whole scheme consistent. Suitable
changes should be made in the case of identical particles. The modified
conditions defining bosons and fermions are given in eqs.\mref{w26}.
It appears that, once the wave function is expressed in terms of 
total and relative variables, they provide the same selection rules as
in
the standard case.


\begin{thebibliography}{99}
\bibitem{b1}
see for example: Takhtajan, L.A. (1992) Lectures on quantum groups, in
{\em
Introduction to Quantum Group and Integrable Massive Models of Quantum
Field Theory},
M.L.Ge\& B.H.Zhao (eds), pp. 69--197, World Scientific, Singapore.
\bibitem{b2}
J.Lukierski, A.Nowicki, H.Ruegg, V.Tolstoy, Phys.Lett. {\bf B264}, 331
(1991);
J.Lukierski, A.Nowicki, H.Ruegg,  Phys.Lett. {\bf B293}, 344 (1993);
S.Giller, P.Kosi\'nski, J.Kunz, M.Majewski, P.Ma\'slanka, Phys.Lett.
{\bf B286}, 57 (1992);
S.Zakrzewski, J.Phys. {\bf A27}, 2075 (1994).
\bibitem{b3}
J.A.De Azcarraga, J.C.Perez Bueno, J.Phys. {\bf A29}, 6353 (1996)
\bibitem{b4}
S.Giller, C.Gonera, P.Kosi\'nski, P.Ma\'slanka, Mod.Phys.Lett. {\bf
A10}, 2757 (1995)
C.Gonera, P.Kosi\'nski, P.Ma\'slanka, M.Tarlini, {\em Projective
representations
of k--Galilei group}, to be published in J.Phys. {\bf A}. 
\bibitem{b5}
J.M.Levi--Leblond, J.Math.Phys. 4 (1963), 776
\bibitem{b6}
J.Lukierski, P.Stichel preprint hep--th 9606170
\end{thebibliography}
\end{document}